\documentclass[10pt,bezier]{article}
\usepackage{amsmath,graphicx,amssymb,amsfonts,color}

\font\bg=cmbx10 scaled\magstep1 \font\Bg=cmbx12 scaled\magstep3
\font\small=cmr8

\newtheorem{newlemma}{{\bf Lemma}}

\newenvironment{lema}{\begin{newlemma}{\hspace{-0.5
em}{\bf.}}}{\end{newlemma}}

\newtheorem{newteorem}{{\bf Theorem}}

\newenvironment{teorem}{\begin{newteorem}{\hspace{-0.5
em}{\bf.}}}{\end{newteorem}}

\newtheorem{newkorolari}{{\bf Corollary}}

\newenvironment{korolari}{\begin{newkorolari}{\hspace{-0.5
em}{\bf.}}}{\end{newkorolari}}

\newtheorem{newdefine}{{\bf Definition}}

\newenvironment{define}{\begin{newdefine}{\hspace{-0.5
em}{\bf.}}}{\end{newdefine}}

\newtheorem{newquestion}{{\bf Question}}

\newtheorem{newkonjek}{{\bf Conjecture}}

\newtheorem{newexample}{{\bf Example}}

\begin{document}
\tolerance=10000 \baselineskip18truept
\newbox\thebox
\global\setbox\thebox=\vbox to 0.2truecm{\hsize
0.15truecm\noindent\hfill}
\def\boxit#1{\vbox{\hrule\hbox{\vrule\kern0pt
     \vbox{\kern0pt#1\kern0pt}\kern0pt\vrule}\hrule}}
\def\qed{\lower0.1cm\hbox{\noindent \boxit{\copy\thebox}}\bigskip}
\def\ss{\smallskip}
\def\ms{\medskip}
\def\bs{\bigskip}
\def\c{\centerline}
\def\nt{\noindent}
\def\ul{\underline}
\def\ol{\overline}
\def\lc{\lceil}
\def\rc{\rceil}
\def\lf{\lfloor}
\def\rf{\rfloor}
\def\ov{\over}
\def\t{\tau}
\def\th{\theta}
\def\k{\kappa}
\def\l{\lambda}
\def\L{\Lambda}
\def\g{\gamma}
\def\d{\delta}
\def\D{\Delta}
\def\e{\epsilon}
\def\lg{\langle}
\def\rg{\rangle}
\def\p{\prime}
\def\sg{\sigma}
\def\ch{\choose}

\newcommand{\ben}{\begin{enumerate}}
\newcommand{\een}{\end{enumerate}}
\newcommand{\bit}{\begin{itemize}}
\newcommand{\eit}{\end{itemize}}
\newcommand{\bea}{\begin{eqnarray*}}
\newcommand{\eea}{\end{eqnarray*}}
\newcommand{\bear}{\begin{eqnarray}}
\newcommand{\eear}{\end{eqnarray}}

\centerline{\Bg  Introduction to  Domination }
 \vspace{.3cm}

\centerline {\Bg Polynomial of a Graph}
\bigskip

\baselineskip12truept \centerline{\bg Saeid
Alikhani$^{a,b,}${}\footnote{\baselineskip12truept\it\small
Corresponding author. E-mail: alikhani@yazduni.ac.ir} and Yee-hock
Peng$^{b,c}$} \baselineskip20truept \centerline{\it $^{a}$Department
of Mathematics} \vskip-8truept \centerline{\it Yazd University }
\vskip-8truept \centerline{\it 89195-741, Yazd, Iran} \vskip-9truept
\centerline{\it $^{b}$Institute for Mathematical Research, and}
\vskip-9truept \centerline{\it $^{c}$Department of Mathematics,}
\vskip-8truept \centerline{\it University Putra Malaysia, 43400 UPM
Serdang, Malaysia} \vskip-8truept \centerline{\it } \vskip-0.2truecm
\nt\rule{12cm}{0.1mm}

\nt{\bg ABSTRACT}
 \medskip

\baselineskip13truept

\noindent{We introduce a domination polynomial of a graph $G$. The
domination polynomial of a graph $G$ of order $n$ is the polynomial
$D(G,x)=\sum_{i=\gamma(G)}^{n} d(G,i) x^{i}$, where $d(G,i)$ is the
number of dominating sets  of $G$ of size $i$, and $\gamma(G)$ is
the domination number of $G$. We obtain some properties  of $D(G,x)$
and its coefficients. Also we compute this polynomial for  some
specific graphs. }

\bs

\nt{\it Keywords:} {\small Domination polynomial; Dominating set;
Unimodal}

\nt{\it Mathematics subject classification:} {\small 05C69, 11B83}

\nt\rule{12cm}{0.1mm}

\baselineskip15truept
\section{Introduction}
Let $G=(V,E)$ be a simple graph. For any vertex $v\in V$, the {\it
open neighborhood} of $v$ is the set $N(v)=\{u \in V|uv\in E\}$ and
the {\it closed neighborhood} of $v$ is the set $N[v]=N(v)\cup
\{v\}$. For a set $S\subseteq V$,
 the open neighborhood of $S$ is $N(S)=\bigcup_{v\in S} N(v)$ and the closed
neighborhood of $S$ is $N[S]=N(S)\cup S$. A set $S\subseteq V$ is a
{\it dominating set} of $G$, if $N[S]=V$, or equivalently, every
vertex in $V\backslash S$ is adjacent to at least one vertex in $S$.
The {\it domination number} $\gamma(G)$ is the minimum cardinality
of a dominating set in $G$. A dominating set with cardinality
$\gamma(G)$ is  called a {\it $\gamma$-set}. For a detailed
treatment of this parameter, the reader is referred
to~\cite{domination}.
 We denote the family of dominating sets of graph $G$ with
cardinality $i$ by ${\cal D}(G,i)$.

\ms

\nt The {\it corona} of two graphs $G_1$ and $G_2$, as defined by
Frucht and Harary in~\cite{harary},
 is the graph
$G=G_1 \circ G_2$ formed from one copy of $G_1$ and $|V(G_1)|$
copies of $G_2$, where the ith vertex of $G_1$ is adjacent to every
vertex in the ith copy of $G_2$. The corona $G\circ K_1$, in
particular, is the graph constructed from a copy of $G$,
 where for each vertex $v\in V(G)$, a new vertex $v'$ and a pendant edge $vv'$ are added.
 The {\it join} of two graphs $G_1$ and $G_2$, denoted by $G_1\vee G_2$,
 is a graph with vertex set  $V(G_1)\cup V(G_2)$
 and edge set $E(G_1)\cup E(G_2)\cup \{uv| u\in V(G_1)$ and $v\in V(G_2)\}$.

\ms

\nt  A finite sequence of real numbers $(a_0,a_1,a_2,\ldots,a_n)$ is
said to be {\it unimodal} if there is some $k\in \{0,1,\ldots,n\}$,
called the {\it mode} of sequence, such that
\[
a_0\leq\ldots\leq a_{k-1}\leq a_k\geq a_{k+1}\geq \ldots\geq a_n;
\]
the mode is {\it unique} if $a_{k-1}<a_k>a_{k+1}$. A polynomial is
called {\it unimodal} if the sequence of its coefficients is
unimodal.

\ms

\nt In the next section, we introduce the domination polynomial and
obtain some of its properties. In Section 3, we study the
coefficients of  the domination polynomials. In the last section, we
investigate the domination polynomial of the graph $G\circ K_1$,
where $G\circ K_1$ is the corona of two graphs $G$ and $K_1$. Also
we show that $D(G\circ K_1,x)$ is unimodal. \ms

\section{Introduction to domination polynomial}

\nt In this section, we state the definition of domination
polynomial and some of its properties.

\begin{define}\label{Definition2}
 Let ${\cal D}(G,i)$ be the family of dominating sets of a graph $G$ with cardinality $i$ and let
$d(G,i)=|{\cal D}(G,i)|$. Then the domination polynomial $D(G,x)$ of
$G$ is defined as
\begin{center}
$D(G,x)=\displaystyle\sum_{i=\gamma(G)}^{|V(G)|} d(G,i) x^{i}$,
\end{center}
where $\gamma(G)$ is the domination number of $G$.
\end{define}

\nt The path $P_4$ on $4$ vertices, for example, has one dominating
set of cardinality $4$, four dominating sets of cardinalities  $3$
and $2$; its domination  polynomial is then $D(P_4,x) =
x^4+4x^3+4x^2$. As another example, it is easy to see that, for
every $n\in \mathbb{N}$,  $D(K_n,x)=(1+x)^n-1$.

\begin{teorem}\label{theorem2.2.2}
If a graph $G$ consists of $m$ components $G_1,\ldots,G_m$, then
$D(G,x)=D(G_1,x)\cdots D(G_m,x)$.
\end{teorem}

\nt{\bf Proof.} It suffices  to prove this theorem for $m=2$. For
$k\geq\gamma(G)$, a dominating set of $k$ vertices in $G$ arises by
choosing a dominating set of $j$ vertices in $G_1$ (for some $j\in
\{\gamma(G_1),\gamma(G_1)+1,\ldots,|V(G_1)|\}$) and a dominating set
of $k-j$ vertices in $G_2$. The number of way of doing this over all
$j=\gamma(G_1),\ldots,|V(G_1)|$ is exactly the coefficient of $x^k$
in $ D(G_1,x)D(G_2,x)$. Hence both side of the above equation have
the same coefficient, so they are identical polynomial.\quad\qed

\nt As a consequence of Theorem~\ref{theorem2.2.2},
 we have the following corollary for the empty graphs:

\begin{korolari}\label{corollary2.2.3}
Let $\overline{K}_n$ be the empty graph with $n$ vertices. Then
$D(\overline{K}_n,x)=x^n$.
\end{korolari}

\nt{\bf Proof.} Since $D(\overline{K}_1,x)=x$, we have the result by
Theorem~\ref{theorem2.2.2}.\quad\qed

\nt Here, we provide a formula for the  domination polynomial of the
join of two graphs.

\begin{teorem}\label{theorem2.2.7}
Let $G_1$ and $G_2$ be  graphs of order  $n_1$ and $n_2$,
respectively. Then
\[
D(G_1\vee
G_2,x)=\Big((1+x)^{n_1}-1\Big)\Big((1+x)^{n_2}-1\Big)+D(G_1,x)+D(G_2,x).
\]
\end{teorem}

\nt{\bf Proof.} Let $i$ be a natural number  $1\leq i \leq n_1+n_2$.
We want to determine $d(G_1\vee G_2,i)$. If $i_1$ and $i_2$ are two
natural numbers such that $i_1+i_2=i$,
 then clearly, for every $D_1\subseteq  V(G_1)$ and $D_2\subseteq V(G_2)$,
such that $|D_j|=i_j$, $j=1,2$, $D_1\cup D_2$ is a dominating set of
$G_1\vee G_2$. Moreover,
 if $D\in  {\cal D}(G_1,i)$,
then $D$ is a dominating set for $G_1\vee G_2$ of size $i$. The same
is true for every $D\in {\cal D}(G_2,i)$. Thus
\[
 D(G_1\vee G_2,x)=\Big((1+x)^{n_1}-1\Big)\Big((1+x)^{n_2}-1\Big)+D(G_1,x)+D(G_2,x).\quad\qed
\]

\nt As a corollary, we have the following formula for the domination
polynomial of the complete bipartite graph $K_{m,n}$, the star
$K_{1,n}$ and the wheel $W_n$.

\begin{korolari}\label{corollary2.2.8}
\begin{enumerate}
\item[(i)]
$D(K_{m,n},x)=((1+x)^m-1)((1+x)^n-1)+x^m+x^n.$
\item[(ii)] $D(K_{1,n},x)=x^n+x(1+x)^n$.
\item[(iii)] If $n\geq 4$, then
$D(W_n,x)=x(1+x)^{n-1}+D(C_{n-1},x)$.
\end{enumerate}

\end{korolari}

\nt {\bf Proof.}
\begin{enumerate}
\item[(i)]
By applying Theorem~\ref{theorem2.2.7} with  $G_1=\overline{K}_n$
and $G_2=\overline{K}_m$, we have the result.
\item[(ii)]
It's suffices to apply Part $(i)$ for $m=1$.
\item[(iii)]
Since for every $n\geq 4$, $W_n=C_{n-1}\vee K_1$, we have the result
by Theorem~\ref{theorem2.2.7}.\quad\qed
\end{enumerate}

\nt In Corollary~\ref{corollary2.2.8}$(iii)$, we have  a
relationship between the domination
 polynomials of wheels and cycles.  For the study of the domination polynomial
 of cycles, the reader is referred to~\cite{saeid}.

\section{Coefficients of domination polynomial}

\nt In this section, we obtain some properties of the coefficients
of the domination polynomial of a graph.

\ms

\nt The following theorem is an easy consequence of the definition
of the domination polynomial.

\begin{teorem}\label{theorem42}
Let $G$ be a graph with $|V(G)|=n$. Then
\begin{enumerate}
\item[(i)] If $G$ is connected, then $d(G,n)=1$ and $d(G,n-1)=n$,
\item[(ii)] $d(G,i)=0$ if and only if $i<\gamma(G)$ or $i>n$.
\item[(iii)] $D(G,x)$ has no constant term.
\item[(iv)] $D(G,x)$ is a strictly increasing function in $[0,\infty)$.
\item[(v)] Let $G$ be a graph and $H$ be any induced subgraph of $G$. Then $deg(D(G,x))\geq deg(D(H,x))$.
\item[(vi)] Zero is a root of $D(G,x)$, with multiplicity $\gamma(G)$.
\end{enumerate}
\end{teorem}

\nt In the following theorem, we want to show that, from the
domination polynomial of a graph $G$, we can obtain the number of
isolated vertices, the number of $K_2$-components and the number of
vertices of degree one in  $G$.

\begin{teorem}\label{theorem2.3.1}
Let $G$ be a graph of order $n$ with $t$ vertices of degree one and
$r$ isolated vertices. If $D(G,x)=\sum_{i=1}^n d(G,i)x^i$
 is its domination polynomial,
then the following hold:

\begin{enumerate}
\item[(i)]
  $r=n-d(G,n-1)$.

\item[(ii)]
 If $G$ has $s$ $K_2$-components,
 then $d(G,n-2)={n\choose 2}-t+s-r(n-1)+{r\choose 2}$.

 \item[(iii)]
 If $G$ has no isolated vertices and $D(G,-2)\neq 0$, then $t={n\choose 2}-d(G,n-2)$.

 \item[(iv)]
 $d(G,1)=\Big|\Big\{v\in V(G)| deg(v)=n-1\Big\}\Big|$.
\end{enumerate}
\end{teorem}

\nt{\bf Proof.}
\begin{enumerate}
\item[(i)]
Suppose that $A\subseteq V(G)$ is the set of  all isolated vertices.
Therefore by assumption,  $|A|=r$. For any vertex $v\in
V(G)\setminus A$, the set  $V(G)\backslash\{v\}$ is a dominating set
of $G$.
 Therefore  $d(G,n-1)=|V(G\backslash A)|=n-r$, and   $r=n-d(G,n-1)$.

\item[(ii)]
Suppose that $D\subseteq V(G)$ is  a  set of  cardinality $n-2$
which is not a dominating set of $G$. We have three cases for $D$:

\nt Case 1. $D=V(G)\backslash\{v,w\}$, where $v$ is an isolated
vertex and $v\in V(G)\backslash \{w\}$. Thus for every isolated
vertex $v$, there are $n-1$ vertices such that
$V(G)\backslash\{v,w\}$ is not a dominating set. Therefore the total
number of $(n-2)$-subsets of $V(G)$ of the form
$V(G)\backslash\{v,w\}$ which is not dominating set ($v$ or $w$ is
an isolated vertex) is $r(n-1)-{r\choose 2}$, since if $v$ and $w$
are isolated vertices, then we count $V(G)\backslash\{v,w\}$ for
both $v$ and $w$.

\nt Case 2. $D=V(G)\backslash\{v,w\}$, for two adjacent vertices $v$
and $w$ with $deg(v)=1$. Since we have $s$ $K_2$-components, the
number of such $\{v, w\}$ is $t-s$ and the proof is complete.

\item[(iii)]
 Since $D(G,-2)\neq 0$, by Theorem~\ref{theorem2.2.2},
 $G$ has no $K_2$-component, and so by Part $(ii)$, we obtain the result.

 \item[(iv)]
 For every $v\in V(G)$, $\{v\}$ is a dominating set if and only if $v$ is adjacent to all vertices.
 The proof is complete.\quad\qed
\end{enumerate}

\ms

\nt We recall that a subset $M$ of $E(G)$ is called a {\it matching}
in $G$ if its elements are not loops and no two of them are adjacent
in $G$; the two ends of an edge in $M$ are said to be {\it matched}
under $M$. A matching $M$ {\it saturates} a vertex $v$, and $v$ is
said to be $M$-{\it saturated} if some edges of $M$ is incident with
$v$; otherwise $v$ is $M$-{\it unsaturated}.

\nt We need the following result to prove
Theorem~\ref{theorem2.2.6}:

\begin{teorem}  {\rm(Hall~\cite{bondy}, p.72)}\label{theorem2.2.5}
Let $G$ be a bipartite graph with bipartition $(X,Y)$. Then $G$
contains a matching that saturates every vertex in $X$ if and only
if for all $S\subseteq X$, $|N(S)|\geq |S|$.
\end{teorem}

\begin{teorem}\label{theorem2.2.6}
Let $G$ be a graph of order $n$. Then for every $0\leq i <
\frac{n}{2}$, we have $d(G,i)\leq d(G,i+1)$.
\end{teorem}
\nt{\bf Proof.} Consider a bipartite graph with two partite sets $X$
and $Y$.
 The vertices of $X$ are dominating sets of $G$ of cardinality $i$, and
the vertices of $Y$ are all $(i+1)$-subsets of $V(G)$. Join a vertex
$A$ of $X$ to a vertex $B$ of $Y$, if $A\subseteq B$. Clearly,  the
degree of each vertex in $X$ is $n-i$. Also for any $B\in Y$, the
degree of $B$ is at most $i+1$. We claim that for any $S\subseteq
X$, $|N(S)|\geq |S|$ and so by Hall's Marriage Theorem, the
bipartite graph has a matching which saturate all vertices of $X$.
By contradiction suppose that there exists $S\subseteq X$ such that
$|N(S)|< |S|$. The number of edges incident with $S$ is $|S|(n-i)$.
Thus by pigeon hole principle, there exists a vertex $B\in Y$ with
degree more than $n-i$. This implies that $i+1\geq n-i+1$. Hence
$i\geq \frac{n}{2}$, a contradiction. Thus for every $S\subseteq X$,
$|N(S)|\geq |S|$ and the claim is proved. Since for every $A\in  X$,
and every $v\in V(G)\backslash A$, $A\cup \{v\}$ is a dominating set
of cardinality $i+1$,  we conclude that $d(G, i+1)\geq d(G,i)$ and
the proof is complete.\quad\qed

\ms

\nt Obviously the result in Theorem~\ref{theorem2.2.6} is useful for
the study of unimodality of domination polynomial. We state the
following conjecture which is similar to the unimodal conjecture for
chromatic polynomial (See~\cite{dong}, p.47): \ms

\nt{\bf Conjecture.} {\it  The domination polynomial of any graph is
unimodal.} \ms

\section{Domination polynomial of $G\circ K_1$}


\nt Let $G$ be any graph with vertex set $\{v_1,\ldots,v_n\}$.
 Add $n$ new vertices $\{u_1,\ldots,u_n\}$ and join $u_i$ to $v_i$ for  $1\leq i \leq n$.
By the definition of the corona of two graphs, we shall  denote
this graph by  $G\circ K_1$. We  study $D(G\circ K_1,x)$ in this
section. Also we show that $D(G\circ K_1,x)$ is unimodal.

\nt We start with the  following lemma:

\begin{lema}\label{lemma5.4.1}
For any graph $G$ of  order $n$, $\gamma(G\circ K_1)=n$.
\end{lema}
\nt{\bf Proof.} If  $D$ is  a dominating set of $G$, then for every
$1\leq i \leq n$, $u_i\in D$ or $v_i\in D$. Therefore $|D|\geq n$.
Since $\{u_1,\ldots,u_n\}$ is a dominating set of $G\circ K_1$, we
have $\gamma(G\circ K_1)=n$.\quad\qed

\nt  By Lemma~\ref{lemma5.4.1}, $d(G\circ K_1,m)=0$ for $m<n$, so we
shall compute $d(G\circ K_1,m)$ for $n\leq m \leq 2n$.

\begin{teorem}\label{theorem5.4.2}
For any graph $G$ of order $n$ and $n\leq m \leq 2n$,
 we have $d(G\circ K_1,m)={n \choose m-n}2^{2n-m}$. Hence $D(G\circ K_1,x)=x^n(x+2)^n$.
\end{teorem}

\nt{\bf Proof.} Suppose that $D$ is a dominating set of $G\circ K_1$
of size $m$. There are ${n\choose m-n}$ possibilities to choose both
vertices of an edge $\{u_i,v_i\}$ for D. Then there remain
$2^{2n-m}$ possibilities to choose the other vertices by selecting
for each pair $\{u_j,v_j\}$ exactly one of these vertices.
 Therefore
\[
d(G\circ K_1,m)={n \choose m-n}2^{2n-m}.\quad\qed
\]

\nt Here, we study the unimodality of the domination polynomial of
$G_n\circ K_1$, where $G_n$ denote a graph with $n$ vertices. Let us
denotes $G\circ K_1$ simply by $G^*$.
  First we state and prove the following theorem for $G_n^*$.

\begin{teorem} \label{theorem5.4.7}
For every $n\in {\mathbb N}$,
$d(G_{3n+2}^*,4n+2)=d(G_{3n+2}^*,4n+3)$.
\end{teorem}
\nt{\bf Proof.} By Theorem~\ref{theorem5.4.2},
$d(G_{3n+2}^*,4n+2)=2^{2n+2}{3n+2\choose n}$
 and $d(G_{3n+2}^*,4n+3)=2^{2n+1}{3n+2\choose n+1}$. Since $2^{2n+2}{3n+2\choose n}=2^{n+1}{3n+2\choose n+1}$,
 we have the result.\quad\qed

\begin{teorem} {\rm(Unimodal Theorem for $G\circ K_1$)} \label{theorem5.4.8}
For every $n\in {\mathbb N}$,
\begin{enumerate}
\item[(i)]
$2^{3n}=d(G_{3n}^*,3n)<d(G_{3n}^*,3n+1)<\ldots<d(G_{3n}^*,4n-1)<d(G_{3n}^*,4n)>d(G_{3n}^*,4n+1)
>\ldots>d(G_{3n}^*,6n-1)>d(G_{3n}^*,6n)=1$,
\item[(ii)]
$2^{3n+1}=d(G_{3n+1}^*,3n+1)<d(G_{3n+1}^*,3n+2)<\ldots<d(G_{3n+1}^*,4n)<d(G_{3n+1}^*,4n+1)>d(G_{3n+1}^*,4n+2)
>\ldots>d(G_{3n+1}^*,6n+1)>d(G_{3n+1}^*,6n+2)=1$,
\item[(iii)]
$2^{3n+2}=d(G_{3n+2}^*,3n+2)<d(G_{3n+2}^*,3n+3)<\ldots<d(G_{3n+2}^*,4n+2)=d(G_{3n+2}^*,4n+3)>d(G_{3n+2}^*,4n+4)>
\ldots>d(G_{3n+2}^*,6n+3)>d(G_{3n+2}^*,6n+4)=1$.
\end{enumerate}
\end{teorem}

\nt{\bf Proof.} Since the proof of all part are similar, we only
prove the part $(i)$:
\begin{enumerate}
\item[(i)]
We shall prove that $d(G_{3n}^*,i)<d(G_{3n}^*,i+1)$ for $3n\leq i
\leq4n-1$ and $d(G_{3n}^*,i)>d(G_{3n}^*,i+1)$ for $4n \leq i\leq
6n-1$. Suppose that $d(G_{3n}^*,i)<d(G_{3n}^*,i+1)$. By
Theorem~\ref{theorem5.4.2} we have
\[
2^{6n-i}{3n\choose i-3n}<2^{6n-i-1}{3n\choose i-3n+1}
\]
So we have $i<4k-\frac{2}{3}$. On the other hand $i\geq 3n$.
Together we have  $3n\leq i \leq4n-1$. Similarly, we have
$d(G_{3n}^*,i)>d(G_{3n}^*,i+1)$ for $4n \leq i\leq 6n-1$. \quad\qed
\end{enumerate}

\nt  By Theorems~\ref{theorem5.4.7} and~\ref{theorem5.4.8}, we
observe that the mode
 for the family $\{D(G_{3n+2}^*,x)\}$ is not unique,
  but for the families $\{D(G_{3n}^*,x)\}$ and $\{D(G_{3n+1}^*,x)\}$, the mode is unique.

\ms

\nt{\bf Remark.}  The unimodality of $D(G^*,x)$
(Theorem~\ref{theorem5.4.8}) also follows
 immediately from the fact that this polynomial has (except zero) only negative real roots. Hence, $D(G^*,x)$ is log-concave and consequently unimodal (see, for example, Wilf~\cite{wilf}).

\ms
\nt {\bf Acknowledgements.} The authors would like to thank the
referee for valuable  comments and suggestions, and  Professor
 Saieed Akbari for his great ideas.

\end{document}